# MONTE CARLO MAXIMUM LIKELIHOOD ESTIMATION FOR DISCRETELY OBSERVED DIFFUSION PROCESSES

By Alexandros Beskos,[1] Omiros Papaspiliopoulos[2]
and Gareth Roberts

*University of Warwick*

This paper introduces a Monte Carlo method for maximum likelihood inference in the context of discretely observed diffusion processes. The method gives unbiased and a.s. continuous estimators of the likelihood function for a family of diffusion models and its performance in numerical examples is computationally efficient. It uses a recently developed technique for the exact simulation of diffusions, and involves no discretization error. We show that, under regularity conditions, the Monte Carlo MLE converges a.s. to the true MLE. For datasize $n \to \infty$, we show that the number of Monte Carlo iterations should be tuned as $\mathcal{O}(n^{1/2})$ and we demonstrate the consistency properties of the Monte Carlo MLE as an estimator of the true parameter value.

**1. Introduction.** We introduce a Monte Carlo method for maximum likelihood inference in the context of discretely observed diffusion processes. The method gives unbiased and a.s. continuous estimates of the likelihood function, which converge uniformly in the parameters to the likelihood function as the Monte Carlo sample size $N$ increases. Additionally, for increasing datasize $n$, the asymptotically optimal algorithm corresponds to selecting $N = \mathcal{O}(n^{1/2})$.

Consider scalar time-homogeneous diffusion processes, defined by stochastic differential equations (SDEs) of the type

(1) $\qquad dV_s = b(V_s; \theta) \, ds + \sigma(V_s; \theta) \, dB_s, \qquad V_s \in \mathcal{V} \subseteq \mathbf{R},$

where $B$ is a Brownian motion. The drift and the diffusion coefficient, $b(\cdot; \theta)$ and $\sigma(\cdot; \theta)$ respectively, are assumed to be known up to a vector of parame-

Received March 2007; revised June 2007.
[1]Supported in part by Greek State Scholarship's Foundation.
[2]Supported in part by EPSRC Grant GR/S61577/01.

ters $\theta \in \Theta \subset \mathbf{R}^d$. Multivariate extensions of this work are addressed in Section 7. The process is assumed to be observed without error at a collection of time instances,

$$\{V_{t_0}, V_{t_1}, \ldots, V_{t_n}\}, \qquad 0 = t_0 < t_1 < \cdots < t_n,$$

and the statistical challenge is to explore characteristics, such as the maximum and the level sets, of the likelihood function

$$\mathfrak{L}_n(\theta) = \prod_{i=1}^{n} L_i(\theta); \qquad L_i(\theta) := p_{\Delta t_i}(V_{t_{i-1}}, V_{t_i}; \theta), \qquad \theta \in \Theta,$$

where $\Delta t_i = t_i - t_{i-1}$ and

(2) $\qquad p_t(v, w; \theta) := \mathrm{P}[V_t \in dw \mid V_0 = v; \theta]/dw, \qquad t > 0, w, v \in \mathcal{V},$

is the transition density of (1). It is well documented (see, e.g., [37] for a recent review) that inference about diffusion models is complicated by the unavailability of the transition density (2), except for limited cases. Thus, inference is carried out using either nonlikelihood approaches, such as estimating equations [10], efficient method of moments [20] and indirect inference [25], or approximate likelihood-based approaches, such as computationally intensive Markov chain Monte Carlo imputation methods [17, 18, 35] and methods which approximate analytically the transition density [1]. Also, approximate Monte Carlo maximum likelihood approaches have been suggested, most notably by [32] and [16].

The method described in this paper, termed the Simultaneous Acceptance Method (SAM), estimates each likelihood contribution $L_i(\cdot)$ independently. To simplify the presentation, let $L(\cdot)$ denote the likelihood contribution of an arbitrary pair of consecutive data points $V_s = v, V_{s+t} = w$, that is, $L(\theta) = p_t(v, w; \theta)$. SAM generates a *random function* $L(\Xi, \theta)$, $\theta \in \Theta$, where $\Xi$ is a random element *independent* of $\theta$, such that for any fixed $\theta \in \Theta$, $\mathbb{E}[L(\Xi, \theta)] = L(\theta)$, and w.p.1, $\theta \mapsto L(\Xi, \theta)$ is continuous. We estimate $L(\cdot)$ by the Monte Carlo functional averages

(3) $\qquad L^N(\cdot) := \frac{1}{N} \sum_{j=1}^{N} L(\Xi^j, \cdot), \qquad \Xi^1, \ldots, \Xi^N$ independent copies of $\Xi$.

Having obtained estimates $L_i^N(\cdot)$ of $L_i(\cdot)$ using (3) independently for each $i = 1, \ldots, n$, we estimate the likelihood function $\mathfrak{L}_n(\cdot)$ by the product

$$\mathfrak{L}_n^N(\cdot) := \prod_{i=1}^{n} L_i^N(\cdot).$$

We demonstrate that our Monte Carlo estimator $\mathfrak{L}_n^N(\cdot)$ is computationally efficient and we detail its theoretical properties. From the computational perspective, the random element consists of standard exponential and Gaussian variables and can be easily simulated, and $L(\Xi, \theta)$ is of a simple calculable form. Moreover, the a.s. continuity of $\mathfrak{L}_n^N(\cdot)$ facilitates the efficient implementation of optimization routines for locating its maximizer. From the theoretical perspective, $\mathfrak{L}_n^N(\theta)$ is an unbiased estimator of $\mathfrak{L}_n(\theta)$ with finite moments for any $\theta \in \Theta$. More importantly, due to the a.s. continuity of $\mathfrak{L}_n^N(\cdot)$, we can resort to the Strong Law of Large Numbers (SLLN) in Banach spaces, to establish that w.p.1 $\mathfrak{L}_n^N(\theta)$ converges to $\mathfrak{L}_n(\theta)$ uniformly in $\theta$, as $N \to \infty$. Uniform convergence to the likelihood implies, under mild conditions, convergence of the maximizers to the MLE, say $\hat{\theta}_n$:

$$(4) \qquad \text{w.p.1}, \qquad \hat{\theta}_n^N := \arg\max_{\theta \in \Theta} \mathfrak{L}_n^N(\theta) \to \hat{\theta}_n \qquad \text{as } N \to \infty.$$

The rate of convergence is shown to be $\mathcal{O}(N^{1/2})$. Also, additional statistics, for example, profile likelihoods and level sets, can be derived as appropriate limits of the corresponding characteristics of the Monte Carlo estimate of the likelihood.

We also investigate the properties of the Monte Carlo method for increasing datasize $n$. The Monte Carlo MLE $\hat{\theta}_n^N$ is a consistent estimator of the true parameter value, say, $\theta_0$, when $N \to \infty$, $n \to \infty$. The optimal algorithm, in terms of computational cost as $n \to \infty$, is obtained when $N = N_n = \mathcal{O}(n^{1/2})$, whence $n^{1/2}(\hat{\theta}_n^N - \theta_0)$ converges in law to a Gaussian random variable.

The construction of the random function $L(\Xi, \cdot)$ is based on a recently developed *retrospective rejection sampling* algorithm called the Exact Algorithm (EA). EA was introduced in [9] and [7]; it returns a draw from any finite-dimensional distribution of the target SDE by rejection sampling with proposals from the Wiener measure. In [8] it was noticed that the transition density of the target diffusion can be written in terms of the acceptance probability of EA, thus, a simple *pointwise* estimator of $L(\theta)$ for any $\theta \in \Theta$ is readily available. The algorithm is termed the Acceptance Method (AM) in [8]. In that paper a simultaneous version of the algorithm is also presented, but for the case of a specific family of diffusions where its development (from pointwise to function estimator) is rather straightforward. In the current paper SAM is applied to a considerably larger family of diffusions and its construction will involve novel couplings of Brownian motion paths. Also, the consistency properties of the method (conditional and unconditional, i.e., for fixed and of increasing number data points respectively) are examined here for the first time.

The applicability of SAM (and AM) is attached to that of EA, whose latest development [7, 8] covers the class of diffusion processes determined by the set of conditions (C$_0$)–(C$_3$) on the drift and diffusion coefficient given in Section 2. Though not universally applicable, SAM can be applied to several diffusion models for which the likelihood function is intractable, for example, to the class of diffusion processes generated by one-to-one transformations of the linear diffusion with multiplicative noise (see, e.g., page 119 of [29]):

$$dV_s = (\theta_1 + \theta_2 V_s)\,ds + (\theta_3 + \theta_4 V_s)\,dB_s, \tag{5}$$

with state space $\mathcal{V} = (\max\{-\theta_3/\theta_4, -\theta_1/\theta_2\}, \infty)$. We will discuss the potential of SAM for more general classes of diffusions. Relevant to this direction are undergoing developments in EA [6].

The paper is organized as follows. In Section 2 we set up the basic notation and provide the random function $L(\Xi, \cdot)$; we state its properties in Theorem 1. In Section 3 we establish a.s. uniform convergence of the estimators $\mathfrak{L}_n^N(\cdot)$ to $\mathfrak{L}_n(\cdot)$, as $N \to \infty$, and discuss some important consequences. In Section 4 we illustrate our method by applying SAM to an example SDE. In Section 5 we allow $n \to \infty$ and present the theory suggesting computational optimality for the choice $N = \mathcal{O}(n^{1/2})$. In Section 6 we prove Theorem 1. In Section 7 we summarize and discuss further ideas. The paper is supplemented with a brief Appendix.

**2. Basic notation and statement of the main result.** In this section we set up the basic notation, construct explicitly the estimator $L(\Xi, \cdot)$ of $L(\cdot) = p_t(v, w; \cdot)$, and state the main, from the applied perspective, result of the paper, which is proven later in Section 6. Throughout the rest of the paper we assume that $\Theta$ is a compact subset of $\mathbf{R}^d$ which contains the unknown MLE $\hat{\theta}_n$. Generally, random variables will be written in capital letters; `typewriter` style will be used to emphasize that the distribution of the variable is independent of $\theta$.

We define

$$\eta(u, \theta) = \int^u \frac{1}{\sigma(z; \theta)}\,dz, \qquad u \in \mathcal{V}, \tag{6}$$

to be any fixed anti-derivative of $1/\sigma(\cdot; \theta)$. For arbitrary $u$, let

$$\alpha(u; \theta) = \frac{b(\eta^{-1}(u, \theta); \theta)}{\sigma(\eta^{-1}(u, \theta); \theta)} - \sigma'(\eta^{-1}(u, \theta); \theta)/2; \qquad A(u, \theta) = \int_0^u \alpha(z; \theta)\,dz,$$

where $\eta^{-1}$ is the inverse of $\eta$. Assuming that $\eta(\cdot, \theta)$ is twice continuously differentiable, Itô's lemma shows that $V_s \mapsto \eta(V_s; \theta) =: X_s$ is the unique transformation (up to a change of sign) that maps $V$ to a process of unit diffusion coefficient; $\alpha(\cdot; \theta)$ is the drift of the transformed process $X$. Throughout the paper we assume that the following conditions hold for all $\theta \in \Theta$:

($C_0$) $\eta(\cdot; \theta)$ is twice continuously differentiable; the law of $X_s = \eta(V_s, \theta)$ has a density w.r.t. the Wiener measure provided by Girsanov's theorem (13);

($C_1$) $\alpha(\cdot; \theta)$ is continuously differentiable;

(C$_2$) $(\alpha^2 + \alpha')(\cdot; \theta)$ is bounded below;
(C$_3$) $(\alpha^2 + \alpha')(\cdot; \theta)$ is bounded above on $(z, \infty)$ for all $z \in \mathbf{R}$.

REMARK 1. The essence of Theorem 1 below is to devise an unbiased estimator of the transition density of the unit diffusion coefficient process $X_s = \eta(V_s, \theta)$ when its drift $\alpha(\cdot; \theta)$ satisfies conditions (C$_1$)–(C$_3$) above. Thus, the case when $(\alpha^2 + \alpha')(\cdot; \theta)$ happens to be bounded above on $(-\infty, z)$ is covered by symmetry. One then needs only to consider instead the (unit diffusion coefficient) process $-\eta(V_s, \theta)$; its drift will then satisfy (C$_1$)–(C$_3$).

We define
$$l(\theta) = \inf_{z \in \mathbf{R}} \frac{(\alpha^2 + \alpha')(z; \theta)}{2}; \qquad \phi(u, \theta) = \frac{(\alpha^2 + \alpha')(u; \theta)}{2} - l(\theta) \geq 0,$$

and $r(u, \theta) = \sup_{z \in (u, \infty)} \phi(z, \theta) < \infty$.

The random element in the estimator $L(\Xi, \theta)$ is $\Xi = (\mathtt{E}, \Psi, \mathtt{Z}, \mathtt{N})$, where $\mathtt{E} \sim \mathrm{Ex}(1)$, $\mathtt{Z} \sim \mathcal{N}(0, 1)$, $\Psi$ is a homogeneous Poisson process on $[0, t]$ with time-ordered points $\mathtt{Y}_j, 1 \leq j \leq \Lambda$, of number $\Lambda \sim \mathrm{Poisson}(\lambda t)$, where the intensity $\lambda$ is specified below, and conditionally on $\Lambda$, the $3 \times \Lambda$-matrix $\mathtt{N} = \{\mathtt{N}_{ij}, 1 \leq i \leq 3, 1 \leq j \leq \Lambda\}$ consists of independent standard Gaussian variables. On the event $\{\Lambda = 0\}$, $\Psi$ and $\mathtt{N}$ contain no elements. The intensity $\lambda$ depends on $\mathtt{E}$ and the observed data, $v, w$, in the following way: we define the functions

(7) $$x = x(\theta) := \eta(v, \theta), \qquad y = y(\theta) := \eta(w, \theta),$$

(8) $$m = m(\mathtt{E}, \theta) = (y + x - \sqrt{2t\mathtt{E} + (y - x)^2})/2,$$

and take $\lambda$ to be any real not less than $\sup_{\theta \in \Theta} r(m, \theta)$. The following theorem gives the random function $L(\Xi, \cdot)$ as a composition of easily computable functions. We define $\mathcal{N}_t(u) := e^{-u^2/(2t)}/\sqrt{2\pi t}$, $u \in \mathbf{R}$, $t > 0$.

THEOREM 1. *Under assumptions* (C$_0$)–(C$_3$) *and the following regularity conditions:*

(Cnt1) *(i)* $\alpha(\cdot; \cdot), \alpha'(\cdot; \cdot)$, *and (ii)* $A(\cdot, \cdot)$ *are continuous on* $\mathbf{R} \times \Theta$,
(Cnt2) *(i)* $\eta(u, \cdot)$, *(ii)* $\eta'(u, \cdot)$ *are continuous for all* $u \in \mathcal{V}$, *(iii)* $l(\cdot)$ *is continuous,*

*the random function* $L(\Xi, \cdot)$ *defined below is such that:*

(i) *for any* $\theta \in \Theta$, $\mathbb{E}[L(\Xi, \theta)] = L(\theta)$,
(ii) *w.p.1* $\theta \mapsto L(\Xi, \theta)$ *is continuous,*
(iii) *there exists a constant* $M$ *such that w.p.1* $|L(\Xi, \theta)| < M$ *for all* $\theta \in \Theta$.

$$L(\Xi, \theta) = |\eta'(w,\theta)|\mathcal{N}_t(y-x)\exp\{A(y,\theta) - A(x,\theta) - l(\theta)t\} \times a(\Xi,\theta);$$

$$a(\Xi, \theta) = \sum_{i=1}^{2} p_i \times \prod_{j=1}^{\Lambda} [1 - \phi(\chi_{ij}, \theta)/\lambda];$$

$$\chi_{ij} = m + \sqrt{\left(\beta_{ij} + \sum_{l=1}^{j} \mathrm{N}_{1l}\gamma_{ijl}\right)^2 + \left(\sum_{l=1}^{j} \mathrm{N}_{2l}\gamma_{ijl}\right)^2 + \left(\sum_{l=1}^{j} \mathrm{N}_{3l}\gamma_{ijl}\right)^2};$$

$$\beta_{ij} = (x-m)\frac{\tau_i - \mathrm{Y}_j}{\tau_i}\mathbb{I}[\mathrm{Y}_j \leq \tau_i] + (y-m)\frac{\mathrm{Y}_j - \tau_i}{t - \tau_i}\mathbb{I}[\mathrm{Y}_j > \tau_i];$$

$$\gamma_{ijl}^2 = \frac{(\tau_i - \mathrm{Y}_j)^2(\mathrm{Y}_l - \mathrm{Y}_{l-1})}{(\tau_i - \mathrm{Y}_l)(\tau_i - \mathrm{Y}_{l-1})}\mathbb{I}[\mathrm{Y}_j \leq \tau_i]$$
$$+ \frac{(t - \mathrm{Y}_j)^2(\mathrm{Y}_l - \mathrm{Y}_{l-1} \vee \tau_i)}{(t - \mathrm{Y}_l)(t - \mathrm{Y}_{l-1} \vee \tau_i)}\mathbb{I}[\mathrm{Y}_l > \tau_i], \qquad \gamma_{ijl} > 0;$$

$$\tau_1 = t\left(1 + \frac{y-m}{x-m}g\right)^{-1}, \qquad \tau_2 = t\left(1 + \frac{y-m}{x-m}g^{-1}\right)^{-1};$$

$$p_1 = \frac{tg\mathrm{E} + 2(x-m)^2}{(1+g)[t\mathrm{E} + 2(x-m)^2]}, \qquad p_2 = 1 - p_1;$$

$$g = 1 + \mathrm{Z}^2/\mathrm{E} - \sqrt{2\mathrm{Z}^2/\mathrm{E} + \mathrm{Z}^4/\mathrm{E}^2} > 0.$$

REMARK 2. Notice that $\chi_{ij}$, $\beta_{ij}$, $\gamma_{ijl}$, $\tau_i$, $p_i$ are all functions of $\theta$, the random element $\Xi$ and the data. Also, $\Xi$ depends on the data points $v, w$ through the Poisson rate $\lambda$.

COROLLARY 1. *Under the assumptions of Theorem 1, the likelihood term $L(\cdot)$ is continuous.*

PROOF. Consider some $\theta \in \Theta$ and a sequence $\{\theta_j\}$ such that $\theta_j \to \theta$. Then, $\lim_{j\to\infty} L(\theta_j) = \lim_{j\to\infty} \mathbb{E}[L(\Xi, \theta_j)] = \mathbb{E}[\lim_{j\to\infty} L(\Xi, \theta_j)] = L(\theta)$, where the first equality holds due to the unbiasedness stated in result (i) of Theorem 1, and the others due to results (ii), (iii) and the bounded convergence theorem. □

**3. Uniform convergence of likelihood estimator.** The Monte Carlo average $L_i^N(\theta)$ in (3) is an unbiased estimator of the likelihood factor $L_i(\theta)$ for any $\theta \in \Theta$, thus, Kolmogorov's Strong Law of Large Numbers (SLLN) implies that w.p.1, $\mathfrak{L}_n^N(\theta) \to \mathfrak{L}_n(\theta)$ as $N \to \infty$. It is known, however, that this pointwise convergence is not strong enough to guarantee convergence of the maximizers $\hat{\theta}_n^N$ of $\mathfrak{L}_n^N(\cdot)$ to the MLE $\hat{\theta}_n$. Similarly, a stronger form of convergence is needed to ensure that several other interesting features of the Monte Carlo functional averages, such as level sets, integrals over subsets and profile likelihoods, will converge to the corresponding features of the likelihood function. A sufficient condition is that the convergence is uniform in $\theta$:

$$(9) \qquad \text{w.p.1,} \qquad \lim_{N \to \infty} \sup_{\theta \in \Theta} |\mathfrak{L}_n^N(\theta) - \mathfrak{L}_n(\theta)| = 0.$$

Questions which lead to essentially equivalent problems to uniform convergence of functional averages occur in theoretical statistics (see, e.g., [11, 38, 39] for the Glivenko–Cantelli problem and [13, 40] for the consistency of the maximum likelihood estimator), dynamical systems and ergodic theory, stochastic optimization, econometrics [2] and Monte Carlo methods [22]; see Chapter 1 of [33] for a review. The general probabilistic framework in which convergence of functional averages is most naturally addressed is *probability on Banach spaces*.

We prove (9) using the following fundamental theorem about SLLN for random elements in an arbitrary separable Banach space. This result first appeared in [30]; see also [5] and [24]. We recall (see, e.g., [24]) that for an arbitrary random element $\mathcal{X}$ in a Banach space $(\mathfrak{X}, \|\cdot\|)$, $\|\mathcal{X}\|$ denotes its norm, and when $\mathbb{E}[\|\mathcal{X}\|] < \infty$, the expectation $\mathbb{E}[\mathcal{X}]$ is defined as the unique element $\mu \in \mathfrak{X}$ such that $T(\mu) = \mathbb{E}[T(\mathcal{X})]$ for every linear function $T : \mathfrak{X} \mapsto \mathbf{R}$ in the dual space of $\mathfrak{X}$.

THEOREM 2. *Let $(\mathfrak{X}, \|\cdot\|)$ be a separable Banach space and $\mathcal{X}$ a random element in $\mathfrak{X}$, such that*

$$\mathbb{E}[\|\mathcal{X}\|] < \infty; \qquad \mathbb{E}[\mathcal{X}] = \mathbf{0}.$$

*If $\mathcal{X}_1, \mathcal{X}_2, \ldots$ are independent copies of $\mathcal{X}$, then*

$$\text{w.p.1,} \qquad \lim_{N \to \infty} \left\| \frac{1}{N} \sum_{j=1}^{N} \mathcal{X}_j \right\| = 0.$$

The uniform convergence in (9) follows easily from the following corollary.

COROLLARY 2. *Let $L(\Xi, \theta)$ be the random function constructed in Theorem 1, and $\Xi^1, \Xi^2, \ldots$, independent copies of $\Xi$. Then*

$$\text{w.p.1,} \qquad \lim_{N \to \infty} \sup_{\theta \in \Theta} \left| \frac{1}{N} \sum_{j=1}^{N} L(\Xi^j, \theta) - L(\theta) \right| = 0.$$

PROOF. Take $(\mathfrak{X}, \|\cdot\|)$ to be the space of continuous real functions on the compact set $\Theta$ equipped with the sup-norm such that, for any element $f \in \mathfrak{X}$, $\|f\| = \sup_{\theta \in \Theta} |f(\theta)|$. This is well known to be a separable Banach space. Theorem 1 has shown that $L(\Xi, \cdot)$ takes values in $\mathfrak{X}$, and $\mathbb{E}[\|L(\Xi, \cdot)\|] < \infty$. Corollary 1 has shown that $L(\cdot) \in \mathfrak{X}$, and by uniqueness of expectation, $\mathbb{E}[L(\Xi, \cdot) - L(\cdot)] = \mathbf{0}$. Applying Theorem 2 to $L(\Xi, \cdot) - L(\cdot)$ yields the result. □

The compactness of $\Theta$ and (9) imply the following result, which validates our Monte Carlo maximum likelihood approach.

COROLLARY 3. *If $\hat{\theta}_n$ is the unique element of $\arg\max_{\theta \in \Theta} \mathfrak{L}_n(\theta)$ and $\{\hat{\theta}_n^N\}_N$ any sequence of maximizers of $\{\mathfrak{L}_n^N(\cdot)\}_N$, then $\lim_{N \to \infty} \hat{\theta}_n^N = \hat{\theta}_n$ w.p.1.*

Note that the same result holds for a sequence of so-called $\epsilon^{(N)}$-maximizers, that is, a sequence $\hat{\theta}_{n,\epsilon}^N$ such that $|\hat{\theta}_n^N - \hat{\theta}_{n,\epsilon}^N| \leq \epsilon^{(N)}$, where $\epsilon^{(N)} \to 0$ as $N \to \infty$. This extension allows for numerically efficient implementations of SAM, as in Section 4.

Certainly, consistency can also be established under the classical approach for showing convergence of the MLE from i.i.d. data (because of the independence of the Monte Carlo samples in the case of SAM) to the true parameter value; see, for example, [19]. It is our impression though that a generalized SLLN provides a natural approach for comprehending strong convergence in enlarged spaces.

Note that a.s. continuity of the random function is not a necessary condition for uniform convergence. For instance, concavity of the functional averages and their expectation would suffice; see [23]. Uniform convergence might hold even when $\Xi^1, \Xi^2, \ldots$, form a stationary (not necessarily ergodic) stochastic process; see, for example, [13, 15, 22, 33]. On the other hand, consistency results as in Corollary 3 could be obtained by establishing epigraphical/hypographical convergence (see, e.g., [4] for relevant results from convex analysis) which is weaker than uniform, although closely related, and requires only semi-continuity properties of $L(\Xi, \cdot)$; see [22] and [13] for investigations in the context of functional averages. In our case, however, concavity does not necessarily hold, but continuity is rather easily established.

Having shown uniform convergence, we can resort to Theorems 5 and 6 of [22] to prove that the Monte Carlo profile log-likelihood and the Monte Carlo level sets converge to the corresponding features of the likelihood function as the number of Monte Carlo samples increases.

3.1. *Asymptotic normality.* Having established strong consistency of $\hat{\theta}_n^N$ as an estimator of the unknown MLE $\hat{\theta}_n$ as $N \to \infty$, we can proceed, under some additional conditions, to prove asymptotic normality for the rescaled sequence $N^{1/2}(\hat{\theta}_n^N - \hat{\theta}_n)$.

We will establish this result by appealing to a general theorem stated below, which concerns the asymptotic normality of the maximizer of an unbiased simultaneous estimator of the likelihood, as the Monte Carlo sample size $N \to \infty$. We first state the result in a general context and then discuss when the conditions it requires are satisfied in our context. We will need the log-likelihood function

$$\ell_n(\theta) := \sum_{i=1}^n \log L_i(\theta)$$

and its estimate

$$\ell_n^N(\theta) := \sum_{i=1}^n \log L_i^N(\theta); \qquad L_i^N(\theta) = \sum_{j=1}^N L_i(\Xi_i^j, \theta)/N.$$

The random elements $\Xi_i^1, \ldots, \Xi_i^N$ are independent copies of, say, $\Xi_i$, which is used for the estimation of the $i$th likelihood term $L_i(\theta) = p_{\Delta t_i}(V_{t_{i-1}}, V_{t_i}; \theta)$. The random elements are independent over the data point index $i = 1, \ldots, n$. We use $\xrightarrow{\mathcal{P}}, \xrightarrow{\mathcal{L}}$ to denote convergence in probability and distribution respectively. All derivatives in the sequel are w.r.t. the parameter argument.

THEOREM 3. *Assume the following:*

(a) $\hat{\theta}_n$ *is unique, in the (assumed nonempty) interior of* $\Theta$.
(b) $\hat{\theta}_n^N \xrightarrow{\mathcal{P}} \hat{\theta}_n$ *as* $N \to \infty$.
(c) $L(\theta) = \mathbb{E}[L(\Xi, \theta)]$ *can be differentiated twice under the expectation sign.*
(d) *There is convergence in distribution*

$$N^{1/2} \nabla \ell_n^N(\hat{\theta}_n) \xrightarrow{\mathcal{L}} \mathcal{N}(0, A_n)$$

*for some covariance matrix $A_n$.*
(e) $B_n = -\nabla^2 \ell_n(\hat{\theta}_n)$ *is positive definite.*
(f) $\nabla^3 \ell_n^N(\theta)$ *is bounded in probability uniformly in a neighborhood of $\hat{\theta}_n$.*

*It is then true that as* $N \to \infty$,

$$-\nabla^2 \ell_n^N(\hat{\theta}_n^N) \xrightarrow{\mathcal{P}} B_n$$

*and*

$$N^{1/2}(\hat{\theta}_n^N - \hat{\theta}_n) \xrightarrow{\mathcal{L}} \mathcal{N}(0, B_n^{-1} A_n B_n^{-1}).$$

PROOF. This is a known result; see, for instance, [19, 22]. Briefly, for the first result note that (b) and the uniform bound in (f) imply, from the mean value theorem, that $\nabla^2 \ell_n^N(\hat{\theta}_n^N) - \nabla^2 \ell_n^N(\hat{\theta}_n) \to 0$ in probability. From (c) and SLLN we get that a.s. $\nabla^2 \ell_n^N(\hat{\theta}_n) \to \nabla^2 \ell_n(\hat{\theta}_n)$. For the second result one takes a first order Taylor expansion with integral remainder for $\nabla \ell_n^N(\hat{\theta}_n^N) - \nabla \ell_n^N(\hat{\theta}_n)$ and rescales by $N^{1/2}$. □

The most restrictive condition of this theorem in the context of SAM is the differentiability of $\theta \mapsto L(\Xi, \theta)$. The formula in Theorem 1 suggests that $L(\Xi, \cdot)$ is typically nondifferentiable at $\{\theta \in \Theta : \tau_i(\theta) = \mathtt{Y}_j, i = 1, 2, j = 1, \ldots, \Lambda\}$. A simple case where this is guaranteed to be an empty set is when the diffusion coefficient $\sigma$ does not depend on $\theta$. If $\sigma(u; \theta) = \sigma(u)$, then also $\eta(u, \theta) = \eta(u)$, and $\theta$ is involved in $L(\Xi, \theta)$ only through the second argument of $\phi(u, \theta)$. Differentiability for $L(\Xi, \cdot)$ can now be implied by straightforward conditions on $\theta \mapsto \alpha(u; \theta)$. Thus, Theorem 3 applies directly to SAM under the assumption that $\sigma(u; \theta) = \sigma(u)$. Note, however, that we have experimentally verified consistency of the same rate $\mathcal{O}(N^{1/2})$ even for the general case; see, for example, the numerical example in the next section. The points of nondifferentiability are of a.s. finite number, and their presence does not seem to effect the result (of local character anyway) of Theorem 3.

We exploit the independence of the transition density estimators $L_i(\Xi_i, \theta)$ over $i = 1, \ldots, n$ to establish condition (d) and identify $A_n$. Notice that

$$N^{1/2} \nabla \ell_n^N(\hat{\theta}_n) = \sum_{i=1}^n \left\{ N^{1/2} \left( \frac{\nabla L_i^N(\hat{\theta}_n)}{L_i^N(\hat{\theta}_n)} - \frac{\nabla L_i(\hat{\theta}_n)}{L_i(\hat{\theta}_n)} \right) \right\}.$$

PROPOSITION 1. *Let $\{X_j, Y_j\}$ be an i.i.d. sequence of vectors $X_j$ and positive scalars $Y_j$ with expected values $\mu_x$ and $\mu_y$ respectively, and running averages $\bar{X}_N = \sum_{j=1}^N X_j/N$, $\bar{Y}_N = \sum_{j=1}^N Y_j/N$. When $N \to \infty$, then*

$$N^{1/2} \left( \frac{\bar{X}_N}{\bar{Y}_N} - \frac{\mu_x}{\mu_y} \right) \xrightarrow{\mathcal{L}} \frac{1}{\mu_y^2} \mathcal{N}(0, \mathrm{Var}(\mu_y X_1 - Y_1 \mu_x)).$$

PROOF. It follows directly from Slutsky's theorems; see [19]. □

Using this proposition, we find that

$$A_n = \sum_{i=1}^n \frac{\mathrm{Var}(L_i(\hat{\theta}_n) \nabla L_i(\Xi_i, \hat{\theta}_n) - L_i(\Xi_i, \hat{\theta}_n) \nabla L_i(\hat{\theta}_n))}{L_i(\hat{\theta}_n)^4}$$

$$= \sum_{i=1}^n \mathrm{Var}\left( \nabla \frac{L_i(\Xi_i, \hat{\theta}_n)}{L_i(\hat{\theta}_n)} \right).$$

$A_n$ can be estimated via the simulated sequences $\{L_i(\Xi_i^j, \hat{\theta}_n^N), \nabla L_i(\Xi_i^j, \hat{\theta}_n^N)\}_j$. Note that $A_n$ increases only linearly with $n$.

**4. Numerical illustration.** We applied SAM to the logistic growth SDE [7]:

$$(10) \qquad dV_s = \delta V_s(1 - c^{-1}V_s)\,ds + \sigma V_s\,dB_s, \qquad V_s \in \mathcal{V} = (0, \infty),$$

which is used to model the evolution of a population in an environment of capacity $c$; $\delta$ is the rate of growth per individual and $\sigma > 0$ a noise parameter. It is known (see, e.g., page 123 of [29]) that (10) is the inverse of the linear SDE with multiplicative noise (5) after making the correspondence $(\theta_1, \theta_2, \theta_3, \theta_4) = (\delta/c, \sigma^2 - \delta, 0, -\sigma)$. Here, we take $\theta = (\delta, c, \sigma)$, and

$$\eta(u, \theta) = -\log(u)/\sigma, \qquad \alpha(u; \theta) = \sigma/2 - \delta/\sigma + \delta/(\sigma c)e^{-\sigma u},$$

$$l(\theta) = \sigma^2/8 - \delta/2, \qquad r(u, \theta) = [(\alpha^2 + \alpha')(u; \theta)/2 - l(\theta)] \vee [\delta^2/(2\sigma^2)].$$

Note that, following Remark 1, $\eta(u, \theta)$ is defined here as the negative of the transformation (6). If $\Theta = [\delta_l, \delta_u] \times [c_l, c_u] \times [\sigma_l, \sigma_u]$, then for any given pair of successive data points $v, w$ of time increment $t$,

$$\lambda = \frac{\delta_u^2}{2\sigma_l^2} \times [(e^{q/2}/c_l - 1)^2 \vee 1]; \qquad q := \log(vw) + \sqrt{2t\sigma_u^2\mathrm{E} + \log^2(w/v)}.$$

It is easy to verify that all conditions of Theorem 1 hold.

In Figure 1 we demonstrate numerically the consistency of the Monte Carlo MLE as an estimator of the unknown MLE. We applied SAM to a dataset of size $n = 1000$ under the specifications $V_0 = 700$, $\theta_0 = (0.1, 1000, 0.1)$ and $\Delta t_i = 1$. (The dataset was simulated using the Exact Algorithm of [7].) We chose $\Theta = [0.03, 0.18] \times [850, 1200] \times [0.09, 0.12]$, outside which a preliminary investigation showed that the likelihood is of negligible value. The maximizers $\hat{\theta}_n^N$ of $\mathfrak{L}_n^N(\cdot)$ for the various values of $N$ in Table 1 were found numerically with the downhill simplex method (see, e.g., Section 10.4 of [34]) up to some precision error $\epsilon^{(N)}$ which decreased with increasing $N$. The initial search point for each maximization was the output of the previous one; for $N = 1$ the initial point was $(0.05, 1150, 0.115)$. Computer implementation in C on a Pentium IV 2.6 GHz processor yielded (in 14 minutes) the sequence of $\epsilon^{(N)}$-maximizers shown in Table 1. The simulated $\Lambda$ was on average, over all consecutive data points and Monte Carlo iterations, 2.004 and never exceeded 11.

In Table 2 we investigate the mean of the asymptotic distribution of the scaled variable $N^{1/2}(\hat{\theta}_n^N - \hat{\theta}_n)$ theoretically demonstrated in Theorem 3. Note that in this context the diffusion coefficient depends on unknown parameters, so $\theta \mapsto L(\Xi, \theta)$ will exhibit points of nondifferentiability with positive probability. Yet, the numerical study suggests agreement with the conclusions of the theorem.

TABLE 1
**Consistency**, $N \to \infty$: The $\epsilon^{(N)}$-maximizers of the estimated likelihood $\mathfrak{L}_n^N(\cdot)$ for a dataset of size $n = 1000$ from the logistic growth model simulated under the parameter values $\theta_0 = (0.1, 1000, 0.1)$. In parenthesis the standard errors found by inverting $-\nabla^2 \ell_n^N(\hat{\theta}_{n,\epsilon}^N)$.

| $N$ | $\delta_{n,\epsilon}^N$ | $c_{n,\epsilon}^N$ | $\sigma_{n,\epsilon}^N$ |
|---|---|---|---|
| 1 | 0.1063 (0.01493) | 1010.0 (30.25) | 0.10051 (0.002354) |
| 2 | 0.1105 (0.01539) | 1009.8 (29.67) | 0.10053 (0.002352) |
| 5 | 0.1115 (0.01562) | 1010.4 (29.31) | 0.10057 (0.002369) |
| 10 | 0.1103 (0.01559) | 1012.3 (29.84) | 0.10060 (0.002364) |
| 20 | 0.1097 (0.01558) | 1012.3 (29.99) | 0.10057 (0.002372) |
| 50 | 0.1100 (0.01563) | 1012.9 (30.02) | 0.10058 (0.002374) |
| 100 | 0.1099 (0.01563) | 1014.2 (30.07) | 0.10058 (0.002371) |
| 200 | 0.1095 (0.01559) | 1014.4 (30.19) | 0.10057 (0.002368) |
| 300 | 0.1096 (0.01560) | 1014.4 (30.19) | 0.10057 (0.002367) |
| Large $N (= 10^4)$ | 0.1096 (0.01561) | 1014.5 (30.18) | 0.10057 (0.002366) |

**5. Unconditional asymptotics, $n \to \infty$.** We find the optimal choice of Monte Carlo iterations $N = N_n$ asymptotically as $n \to \infty$. In accordance, we will now write $\hat{\theta}_n^{N_n}$ for the sequence of Monte Carlo maximizers (4) and $L_i^{N_n}(\theta)$ for the unbiased estimators of $L_i(\theta)$. The data points are now treated as random variables. The analysis follows a rather intuitive approach. The objective is to provide a rule for the selection of $N$ but, at the same time, avoid the mathematical rigor that would take up the space of a full paper; see, for example, [27] or [12]. So, we will state the main result in Theorem 4 under reasonable ergodicity assumptions on the diffusion dynamics without insisting on technical details. Note that the precise formula for the estimator $L(\Xi, \theta)$ of $L(\theta)$ is not essential for this section: given that Assumptions 1–4 below are satisfied, Theorem 4 holds for any given *unbiased* estimator of $L(\theta)$. As in Section 3.1, the asymptotic normality result in Theorem 4 requires differentiability of $\theta \mapsto L(\Xi, \theta)$ and will apply directly to SAM under the known diffusion coefficient assumption, $\sigma(u; \theta) = \sigma(u)$.

TABLE 2
**Asymptotic Unbiasedness**, $N \to \infty$: The sample means from 1000 realizations of $N^{1/2}(\hat{\theta}_n^N - \hat{\theta}_n)$ for various $N$. In parenthesis the corresponding standard errors. The data were of size $n = 250$; $\theta_0 = (0.1, 1000, 0.1)$. The MLE $\hat{\theta}_n$ was found using large $N$ $(= 10^4)$.

|  | $N^{1/2}(\hat{\delta}_n^N - \hat{\delta}_n)$ | $N^{1/2}(\hat{c}_n^N - \hat{c}_n)$ | $N^{1/2}(\hat{\sigma}_n^N - \hat{\sigma}_n)$ |
|---|---|---|---|
| $N = 25$ | $-129\text{e}{-5}$ ($30\text{e}{-5}$) | 0.825 (0.320) | $-2.28\text{e}{-5}$ ($0.79\text{e}{-5}$) |
| $N = 50$ | $-81\text{e}{-5}$ ($30\text{e}{-5}$) | 0.414 (0.321) | $-0.47\text{e}{-5}$ ($0.78\text{e}{-5}$) |
| $N = 100$ | $-46\text{e}{-5}$ ($30\text{e}{-5}$) | 0.092 (0.317) | $1.40\text{e}{-5}$ ($0.81\text{e}{-5}$) |

We assume that the diffusion is ergodic, with $\pi$ denoting its invariant density corresponding to the correct parameter value $\theta_0$, and that the data are equidistant, that is, $t_i = i\Delta$ for some $\Delta > 0$. We define the family of mappings which satisfy a Law of Large Numbers criterion:

$$\Im = \left\{ f \mid f : \mathbf{R}^2 \to \mathbf{R}, \sum_{i=1}^n f(V_{(i-1)\Delta}, V_{i\Delta})/n \xrightarrow{\mathcal{P}} \mathbb{E}_\pi[f(V_0, V_\Delta)] \right\},$$

where $\mathbb{E}_\pi[\cdot]$ is expectation in stationarity, that is, $(V_0, V_\Delta) \sim \pi(dx)p_\Delta(x, dy; \theta_0)$.

We follow the approach of [27]. In that paper the discretization increment of a diffusion approximation is adjusted to the datasize $n \to \infty$. In our case, had it been possible to construct an unbiased estimator of the *log*-transition density, then the Monte Carlo error would be averaged out for increasing $n$ and $\sqrt{n}$-consistency of $\hat{\theta}_n^{N_n}$ (as an estimator of $\theta_0$) would follow even for fixed $N$; see Theorem 2 of [27]. We will now allow $N = N_n \to \infty$ as $n \to \infty$ and identify the magnitude of the log-bias through an error expansion. For simplicity, we set

$$\zeta_i = \zeta_{i,N} := \nabla \log L_i^N(\theta_0) - \nabla \log L_i(\theta_0),$$

and $\mathcal{F}_i = \sigma(V_{(i-1)\Delta}, V_{i\Delta})$, for $i = 1, \ldots, n$. Recall that $d$ is the dimensionality of the parameter vector.

ASSUMPTION 1. *There exist* $\psi : \mathbf{R}^2 \mapsto \mathbf{R}^d$, $g : \mathbf{R}^2 \mapsto \mathbf{R}^{d \times d}$, $h : \mathbf{R}^2 \mapsto \mathbf{R}$ *with scalar components in* $\Im$ *such that*

$$\left| \mathbb{E}[\zeta_{i,N} \mid \mathcal{F}_i] - \psi(V_{(i-1)\Delta}, V_{i\Delta}) \frac{1}{N} \right| + \left| \mathbb{E}[\zeta_{i,N} \zeta_{i,N}^\top \mid \mathcal{F}_i] - g(V_{(i-1)\Delta}, V_{i\Delta}) \frac{1}{N} \right|$$
$$\leq h(V_{(i-1)\Delta}, V_{i\Delta}) o(1/N).$$

The assumption is not on the expansion itself, but on the regularity of the coefficients $\psi, g, h$. Indeed, for given $V_{(i-1)\Delta}, V_{i\Delta}$, the $\mathcal{O}(1/N)$-bias follows from a second-order Taylor expansion for $\zeta_{i,N}$. Analytically, for sequence $\{X_j, Y_j\}$ as in Proposition 1, one can write

$$\frac{\bar{X}_N}{\bar{Y}_N} - \frac{\mu_x}{\mu_y} = \frac{1}{\mu_y}(\bar{X}_N - \mu_x) - \frac{\mu_x}{\mu_y^2}(\bar{Y}_N - \mu_y) - \frac{1}{\mu_y^2}(\bar{X}_N - \mu_x)(\bar{Y}_N - \mu_y)$$
$$+ \frac{\mu_x}{\mu_y^3}(\bar{Y}_N - \mu_y)^2 + R(\bar{X}_N, \bar{Y}_N, \mu_x, \mu_y),$$

for some lower order residual $R(\bar{X}_N, \bar{Y}_N, \mu_x, \mu_y)$. Making the correspondence $\bar{X}_N \leftrightarrow \nabla L_i^N(\theta_0)$, $\bar{Y}_N \leftrightarrow L_i^N(\theta_0)$, assuming that $\mathbb{E}[L_i(\Xi_i, \theta)] = L_i(\theta)$ can be differentiated under the integral sign, and taking expectations conditionally on the data, $\psi(V_{(i-1)\Delta}, V_{i\Delta})$ in Assumption 1 is analytically identified as

$$-\frac{1}{L_i^2(\theta_0)} \operatorname{Cov}(\nabla L_i(\Xi_i, \theta_0), L_i(\Xi_i, \theta_0)) + \frac{\nabla L_i(\theta_0)}{L_i^3(\theta_0)} \operatorname{Var}(L_i(\Xi_i, \theta_0)).$$

The corresponding expressions we can obtain for $g$, $h$ are much more complicated.

ASSUMPTION 2. *The matrix $A(\theta) := \mathbb{E}_\pi[-\nabla^2 \log p_\Delta(V_0, V_\Delta; \theta)]$ is positive definite for any $\theta \in \Theta$ and*

$$\frac{-\nabla^2 \ell_n^{N_n}(\theta)}{n} \xrightarrow{\mathcal{P}} A(\theta)$$

*uniformly for $\theta$ in a neighborhood of $\theta_0$.*

Under the standard maximum likelihood assumption for convergence in probability of $\nabla^2 \ell_n(\theta)/n$ to $A(\theta)$ uniformly in a neighborhood of $\theta_0$, it remains to explain the same mode of convergence for $\{\nabla^2 \ell_n^{N_n}(\theta) - \nabla^2 \ell_n(\theta)\}/n$ toward 0. We define

$$\dot{\zeta}_i = \dot{\zeta}_{i,N} = \nabla^2 \log L_i^N(\theta) - \nabla^2 \log L_i(\theta).$$

From Lemma 9 of [21], used also in [27], the following conditions imply the required convergence for fixed $\theta$:

(11) $$\frac{\sum_{i=1}^n \mathbb{E}[\dot{\zeta}_{i,N_n} \mid \mathcal{F}_i]}{n} \xrightarrow{\mathcal{P}} 0; \qquad \frac{\sum_{i=1}^n \mathbb{E}[\dot{\zeta}_{i,N_n} \dot{\zeta}_{i,N_n}^\top \mid \mathcal{F}_i]}{n^2} \xrightarrow{\mathcal{P}} 0.$$

Similarly to Assumption 1, a Taylor expansion can provide estimates:

$$|\mathbb{E}[\dot{\zeta}_{i,N} \mid \mathcal{F}_i]| + |\mathbb{E}[\dot{\zeta}_{i,N} \dot{\zeta}_{i,N}^\top \mid \mathcal{F}_i]| \leq \dot{R}_\theta(V_{(i-1)\Delta}, V_{i\Delta}) \mathcal{O}(1/N).$$

Assuming that $\dot{R}_\theta \in \Im$, and since $N_n \to \infty$, these estimates imply (11). It is technically much harder to illustrate convergence uniformly in a neighborhood of $\theta_0$. Following Theorem 2 in [27], it suffices (assuming stationarity) to obtain a bound on the $L^{2d+1}$-norm of $\dot{\zeta}_{i,N}(\theta)$ and a Lipschitz condition for $\theta \mapsto \dot{\zeta}_{i,N}(\theta)$ uniformly in $N$, again in the $L^{2d+1}$-norm. We avoid further details.

ASSUMPTION 3. *The following two sequences converge in probability to 0:*

$$\frac{\sum_{i=1}^n \zeta_{i,N_n}^2 - \mathbb{E}[\zeta_{i,N_n}^2 \mid \mathcal{F}_i]}{n}, \qquad \frac{\sum_{i=1}^n \zeta_{i,N_n} \mathbb{E}[\zeta_{i,N_n} \mid \mathcal{F}_i] - \mathbb{E}^2[\zeta_{i,N_n} \mid \mathcal{F}_i]}{n}.$$

Note that the summands have zero expectation. This can therefore be interpreted as an ergodicity assumption on the diffusion dynamics upon the consideration of the enlarged $\sigma$-algebra from the Monte Carlo scheme.

ASSUMPTION 4. *The following weak convergence holds:*
$$\frac{\nabla \ell_n(\theta_0)}{\sqrt{n}} \xrightarrow{\mathcal{L}} \mathcal{N}(0, V),$$
*where*
$$V = \mathbb{E}_\pi[\nabla \log p_\Delta(V_0, V_\Delta; \theta_0) \nabla \log p_\Delta(V_0, V_\Delta; \theta_0)^\top].$$

This is a standard ergodicity assumption for maximum likelihood inference for discretely observed diffusions. Of course, $V = A(\theta_0)$, assuming exchangeability of differentiation and integration at $\mathbb{E}_\pi[\nabla^2 p_\Delta(V_0, V_\Delta; \theta_0)]$.

THEOREM 4. *Let $N = N_n$ grow to infinity with $n \to \infty$. Then*
$$\hat{\theta}_n^{N_n} \xrightarrow{\mathcal{P}} \theta_0.$$
*Also, if $\theta_0$ is in the interior of $\Theta$ and Assumptions 1–4 hold, then:*

(i) *if $\lim_{n \to \infty} \sqrt{n}/N_n = 0$ then $\sqrt{n}(\hat{\theta}_n^{N_n} - \theta_0) \xrightarrow{\mathcal{L}} \mathcal{N}(0, A^{-1})$,*
(ii) *if $\lim_{n \to \infty} \sqrt{n}/N_n = c \in (0, \infty)$ then $\sqrt{n}(\hat{\theta}_n^{N_n} - \theta_0) \xrightarrow{\mathcal{L}} \mathcal{N}(c\mu, A^{-1})$,*
(iii) *if $\lim_{n \to \infty} \sqrt{n}/N_n = \infty$ then $N_n(\hat{\theta}_n^{N_n} - \theta_0) \xrightarrow{\mathcal{P}} \mu$,*

*for $A = A(\theta_0)$ and $\mu = A^{-1}\mathbb{E}_\pi[\psi(V_0, V_\Delta)]$.*

PROOF. Consistency of the Monte Carlo MLE can be proved in line with [12]. We proceed directly at the proof of the asymptotic results for the various scalings of $N_n$. We use a Taylor expansion as in the classical MLE theory. The additional complexity is due to the presence of the Monte Carlo scheme. The basic equation is the following:

$$(12) \quad \nabla \ell_n^N(\hat{\theta}_n^N) = \nabla \ell_n^N(\theta_0) + \int_0^1 \nabla^2 \ell_n^N(\theta_0 + \rho(\hat{\theta}_n^N - \theta_0)) \, d\rho \times (\hat{\theta}_n^N - \theta_0),$$

where the term $\nabla \ell_n^N(\hat{\theta}_n^N) = 0$ can be ignored. From the consistency of $\hat{\theta}_n^{N_n}$ and Assumption 2,

$$\frac{-\int_0^1 \nabla^2 \ell_n^{N_n}(\theta_0 + \rho(\hat{\theta}_n^{N_n} - \theta_0)) \, d\rho}{n} \xrightarrow{\mathcal{P}} A.$$

Consider now the remaining gradient term. We rewrite
$$\nabla \ell_n^N(\theta_0) = \{\nabla \ell_n^N(\theta_0) - \nabla \ell_n(\theta_0)\} + \nabla \ell_n(\theta_0).$$

Assumption 4 controls $\nabla \ell_n(\theta_0)$. For the last term $\{\nabla \ell_n^N(\theta_0) - \nabla \ell_n(\theta_0)\}$ we employ a martingale decomposition as in [27]. We rewrite

$$\nabla \ell_n^N(\theta_0) - \nabla \ell_n(\theta_0) = \sum_{i=1}^n \{\zeta_{i,N} - \mathbb{E}[\zeta_{i,N} \mid \mathcal{F}_i]\} + \sum_{i=1}^n \mathbb{E}[\zeta_{i,N} \mid \mathcal{F}_i].$$

Assumption 1 controls the extreme right term. Under Assumptions 1 and 3,

$$\frac{\sum_{i=1}^{n}(\zeta_{i,N_n} - \mathbb{E}[\zeta_{i,N_n} \mid \mathcal{F}_i])^2}{n} \xrightarrow{\mathcal{P}} 0,$$

so the martingale CLT in Theorem 3.2 of [26] implies

$$\frac{\sum_{i=1}^{n}\{\zeta_{i,N_n} - \mathbb{E}[\zeta_{i,N_n} \mid \mathcal{F}_i]\}}{\sqrt{n}} \xrightarrow{\mathcal{L}} 0.$$

For (i) and (ii), multiply both sides of (12) with $n^{-1/2}$, whence, from Assumption 1, the log-bias term $n^{-1/2}\sum_{i=1}^{n}\mathbb{E}[\zeta_{i,N_n} \mid \mathcal{F}_i]$ will converge in probability either to 0 if $n^{1/2}/N_n \to 0$ or to $c\mu$ if $n^{1/2}/N_n \to c$. For (iii), multiply both sides of (12) with $N_n/n$. □

The theorem shows that, as long as $N = o(n^{1/2})$, the Monte Carlo MLE has the same asymptotic behavior with the unknown MLE.

**6. Proof of Theorem 1.** The notation in this section follows the definitions in Section 2. We prove the result in several stages.

6.1. *Diffusion transformation and densities.* We consider the modified process $X_s = \eta(V_s, \theta)$. Condition (C$_0$) allows the application of Itô's lemma to show that $X_s$ is the solution of the unit diffusion coefficient SDE:

$$dX_s = \alpha(X_s; \theta) \, ds + dB_s.$$

Let $\tilde{p}_t(\cdot, \cdot; \theta)$ be the transition density of $X$ defined analogously to (2). Recall that we have defined $x = x(\theta) = \eta(v, \theta)$ and $y = y(\theta) = \eta(w, \theta)$. A standard change-of-variables argument yields

$$L(\theta) \equiv p_t(v, w; \theta) = |\eta'(w, \theta)| \cdot \tilde{p}_t(x, y; \theta).$$

Let $\mathbb{Q}_\theta$ be the law of the paths of $X$ on $[0, t]$ conditioned to begin at $X_0 = x$ and finish at $X_t = y$. We call such a path "a bridge from $(0, x)$ to $(t, y)$." Let $\mathbb{W}_\theta$ be the distribution of the Brownian bridges (BBs) from $(0, x)$ to $(t, y)$; a random process distributed according to $\mathbb{W}_\theta$ will be denoted by $W$. Both $\mathbb{Q}_\theta$ and $\mathbb{W}_\theta$ apply on the space $C$ of continuous mappings from $[0, t]$ to **R** equipped with the corresponding cylinder $\sigma$-algebra; we denote by $\omega$ a typical element of $C$. We have assumed in (C$_0$) absolute continuity of the law of (unconditional) paths of $X$ w.r.t. the Wiener measure with density given by Girsanov's theorem (see, e.g., [31]):

$$\text{(13)} \qquad \exp\biggl\{\int_0^t \alpha(\omega_s; \theta) \, d\omega_s - \tfrac{1}{2}\int_0^t \alpha^2(\omega_s; \theta) \, ds\biggr\}.$$

Bayes' theorem and an application of Itô's formula gives the following expression for the corresponding bridge density:

$$\frac{d\mathbb{Q}_\theta}{d\mathbb{W}_\theta}(\omega) = \frac{\mathcal{N}_t(y-x)}{\tilde{p}_t(x,y;\theta)} \exp\bigg\{A(y,\theta) - A(x,\theta) - \int_0^t \frac{1}{2}(\alpha^2 + \alpha')(\omega_s;\theta)\,ds\bigg\}.$$

Taking expectations at both sides w.r.t. $\mathbb{W}_\theta$, using ($C_2$) and re-arranging, we get that

$$\tilde{p}_t(x,y;\theta) = \mathcal{N}_t(y-x)\exp\{A(y,\theta) - A(x,\theta) - l(\theta)t\} \times a(\theta),$$

where

$$a(\theta) = \mathbb{E}_{\mathbb{W}_\theta}\bigg[\exp\bigg\{-\int_0^t \phi(\omega_s,\theta)\,ds\bigg\}\bigg] \leq 1.$$

A heuristic description of the remainder of the proof is as follows. We initially derive an unbiased estimator $a(\Xi,\theta)$ of $a(\theta)$, where the distribution of the random element $\Xi$ depends on $\theta$. We then construct jointly the family $\{\Xi; \theta \in \Theta\}$ by expressing $\Xi = f(\Xi, \mathtt{X}, \theta)$ for $\Xi$ given in Theorem 1, $\mathtt{X}$ some other random element also independent of $\theta$ and $f$ an appropriately specified function. Finally, we show that $a(\Xi,\theta)$ defined in Theorem 1 is given as $a(\Xi,\theta) = \mathbb{E}[a(f(\Xi,\mathtt{X},\theta),\theta)\,|\,\Xi]$, where $\mathtt{X}$ is integrated out to ensure a.s. continuity of the mapping $\theta \mapsto a(\Xi,\theta)$.

6.2. *Connection with exact diffusion bridge simulation.* In [8] it is noticed that $a(\theta)$ coincides with the acceptance probability of a rejection sampling algorithm for the simulation of paths with law $\mathbb{Q}_\theta$. The algorithm, developed in [7] and termed the Exact Algorithm (EA), proposes paths from $\mathbb{W}_\theta$ and accepts them according to the density ratio

(14) $$\frac{d\mathbb{Q}_\theta}{d\mathbb{W}_\theta}(\omega) \propto \exp\bigg\{-\int_0^t \phi(\omega_s,\theta)\,ds\bigg\} \leq 1.$$

Let $m = \inf\{\omega_s, s \in [0,t]\}$ be the minimum of $\omega$. Condition ($C_3$) implies that $\phi(\omega_s,\theta)/r(m,\theta) \leq 1$, for all $s \in [0,t]$. When $\omega$ is a realization of the Brownian bridge $W \sim \mathbb{W}_\theta$, the distribution of its minimum is given in terms of the Rayleigh distribution (see, e.g., [36]) and can be simulated precisely according to (8) using the exponential random variable $\mathtt{E}$. Notice that if $\Phi$ is a Poisson process of intensity $r(m,\theta)$ on $[0,t] \times [0,1]$ and $N$ the number of the points of $\Phi$ lying below the curve $s \mapsto \phi(\omega_s,\theta)/r(m,\theta)$, then

$$\mathrm{P}[N=0\,|\,\omega] = \exp\bigg\{-\int_0^t \phi(\omega_s,\theta)\,ds\bigg\}.$$

Thus, in [7] the following simulation algorithm is suggested. Let $\Psi$ be the projection of $\Phi$ on the time-axis with time-ordered points $Y_j, 1 \leq j \leq \Lambda$, where $\Lambda \sim \mathrm{Poisson}(r(m,\theta)t)$.

The Exact Algorithm (EA).

1. Simulate $\mathtt{E} \sim \text{Ex}(1)$, set $m = (y + x - \sqrt{2t\mathtt{E} + (y-x)^2})/2$.
2. Simulate $\Psi$.
3. Simulate $\{W_{Y_j}, 1 \leq j \leq \Lambda\}$, so that its minimum is $m$.
4. With probability $1 - \prod_{j=1}^{\Lambda}[1 - \phi(W_{Y_j}, \theta)/r(m, \theta)]$ goto Step 1.
5. Output all information about $W$.

Step 3 requires the simulation of the skeleton $\{W_{Y_j}, 1 \leq j \leq \Lambda\}$ with a given minimum. This can be achieved (see Section 6.4 below) using existing theory on the decomposition of the Brownian path at its minimum; see, for example, Proposition 2 of [3]. Therefore, a *pointwise* unbiased estimator of $a(\theta)$ is readily available:

$$(15) \qquad a(\Xi, \theta) = \prod_{j=1}^{\Lambda}[1 - \phi(W_{Y_j}, \theta)/r(m, \theta)],$$

where $\Xi = (\mathtt{E}, \Psi, \{W_{Y_j}, 1 \leq j \leq \Lambda\})$. This is the Acceptance Method proposed in [8].

6.3. *Coupling of the Poisson processes.* We define the joint structure of the collection of random variables $\{\Psi; \theta \in \Theta\}$ using the thinning property of the Poisson process (see, e.g., Section 5 of [28]). Let $\bar{\Psi}$ be a Poisson process of rate $\lambda \geq \sup_{\theta \in \Theta} r(m(\mathtt{E}, \theta), \theta)$ on the interval $[0, t]$ with time-ordered points $\mathtt{Y}_j, 1 \leq j \leq \bar{\Lambda}$ and $\mathtt{U} = (\mathtt{U}_1, \ldots, \mathtt{U}_{\bar{\Lambda}})$ a vector of i.i.d. variables, $\mathtt{U}_1 \sim \text{Un}[0, 1]$. We set $\Psi = \{\mathtt{Y}_j \in \bar{\Psi} : \mathtt{U}_j < r(m, \theta)/\lambda, 1 \leq j \leq \bar{\Lambda}\}$. Then the right-hand side of (15) rewrites as

$$\prod_{j=1}^{\bar{\Lambda}} \{1 - \mathbb{I}[\mathtt{U}_j < r(m, \theta)/\lambda] \cdot \phi(W_{\mathtt{Y}_j}, \theta)/r(m, \theta)\}.$$

After integrating out $\mathtt{U}$, we have the following *pointwise* unbiased estimator of $a(\theta)$:

$$(16) \qquad a(\Xi, \theta) = \prod_{j=1}^{\bar{\Lambda}}[1 - \phi(W_{\mathtt{Y}_j}, \theta)/\lambda],$$

where now $\Xi = (\mathtt{E}, \bar{\Psi}, \{W_{\mathtt{Y}_j}, 1 \leq j \leq \bar{\Lambda}\})$.

6.4. *Coupling of the proposed paths.* For a path $\omega$ with given minimum $m$, let $\tau = \{s \in [0, t] : \omega_s = m\}$ be the instance when $m$ is attained. As shown in the Appendix, when $\omega$ is a realization of $W$, the time instance of the minimum can be simulated conditionally on $m$ as follows:

$$(17) \qquad \tau = \mathbb{I}[\mathtt{V} \leq p_1]\tau_1 + \mathbb{I}[\mathtt{V} > p_1]\tau_2,$$

for $\mathtt{V} \sim \mathrm{Un}[0,1]$, where we recall that $p_1, \tau_1, \tau_2$ are given in terms of $\mathtt{E}, \mathtt{Z}$ and $\theta$. A Brownian bridge $W$ with a given minimum $m$ attained at time $\tau$ can be constructed in terms of two Bessel bridges which can in turn be obtained through six independent standard [from $(0,0)$ to $(1,0)$] BBs; see [3]. In detail, the theory allows us to construct $W$ at the time instances $\mathtt{Y}_j, 1 \leq j \leq \Lambda$, in the following way:

$$(18) \quad W_{\mathtt{Y}_j} = m + \begin{cases} \sqrt{\tau}\sqrt{\left(\dfrac{(x-m)(\tau - \mathtt{Y}_j)}{\tau^{3/2}} + \mathtt{W}_{1,\mathtt{Y}_j/\tau}\right)^2 + \mathtt{W}_{2,\mathtt{Y}_j/\tau}^2 + \mathtt{W}_{3,\mathtt{Y}_j/\tau}^2}, \\ \hfill \text{if } \mathtt{Y}_j \leq \tau, \\ \sqrt{t-\tau}\Bigg(\left(\dfrac{(y-m)(\mathtt{Y}_j - \tau)}{(t-\tau)^{3/2}} + \mathtt{W}_{4,(\mathtt{Y}_j-\tau)/(t-\tau)}\right)^2 \\ \hfill + \mathtt{W}_{5,(\mathtt{Y}_j-\tau)/(t-\tau)}^2 + \mathtt{W}_{6,(\mathtt{Y}_j-\tau)/(t-\tau)^2}\Bigg)^{1/2}, \\ \hfill \text{if } \mathtt{Y}_j > \tau, \end{cases}$$

where $\mathtt{W} = \{(\mathtt{W}_{i,s}, 0 \leq s \leq 1), 1 \leq i \leq 6\}$ is the required collection of standard BBs. We define $W_{i,\mathtt{Y}_j}$ as $W_{\mathtt{Y}_j}$ in (18) under the specification $\tau = \tau_i$, for $i = 1, 2$ and $1 \leq j \leq \Lambda$. Under this convention, substituting (17), (18) into (16) and integrating out $\mathtt{V}$, we rewrite the function we wish to estimate as

$$a(\theta) = \mathbb{E}\left\{\mathbb{E}\left[\prod_{j=1}^{\Lambda}[1 - \phi(W_{\mathtt{Y}_j}, \theta)/\lambda] \mid \mathtt{E}, \mathtt{Z}, \mathtt{W}, \Psi\right]\right\}$$

$$= \sum_{i=1}^{2} \mathbb{E}\left\{\prod_{j=1}^{\Lambda}[1 - \phi(W_{i,\mathtt{Y}_j}, \theta)/\lambda] \times p_i\right\} = a_1(\theta) + a_2(\theta),$$

where $a_1(\theta), a_2(\theta)$ are defined in the obvious way. In the sequel, we construct estimators $a_i(\Xi, \theta)$ of $a_i(\theta), i = 1, 2$, and take $a(\Xi, \theta) = a_1(\Xi, \theta) + a_2(\Xi, \theta)$.

For each $\mathtt{Y}_j$, we need to simulate the bridges $\mathtt{W}_1, \mathtt{W}_2, \mathtt{W}_3$ at the instance $\mathtt{Y}_j/\tau_i$ or simulate $\mathtt{W}_4, \mathtt{W}_5, \mathtt{W}_6$ at the instance $(\mathtt{Y}_j - \tau_i)/(t - \tau_i)$ if $\mathtt{Y}_j \leq \tau_i$ or $\mathtt{Y}_j > \tau_i$ respectively, $1 \leq j \leq \Lambda$; see also Figure 1 for a graphical illustration. Thus, in total we need $3 \times \Lambda$ simulations, which can be done using the same number of Gaussian random variables as we now show. We can simulate a standard BB $(\mathtt{W}_s; 0 \leq s \leq 1)$ at some time instances $0 < s_1 < \cdots < s_d < 1$ using the formula (derived using standard properties of BB)

$$(19) \qquad \mathtt{W}_{s_j} = \sum_{l=1}^{j} N_l \sqrt{\dfrac{(1-s_j)^2(s_l - s_{l-1})}{(1-s_{l-1})(1-s_l)}}, \qquad 1 \leq j \leq d,$$

where $N_1, N_2, \ldots, N_d$ are i.i.d. with $N_1 \sim \mathcal{N}(0,1)$. Consider now the collection $\mathtt{N} = \{\mathtt{N}_{kj}, 1 \leq k \leq 3, 1 \leq j \leq \Lambda\}$ of independent $\mathcal{N}(0,1)$ variables. We will use $\{N_{kj}, 1 \leq j \leq \Lambda\}$ for the generation of the required instances of the bridges $W_k$ and $W_{k+3}$, $1 \leq k \leq 3$; see Figure 1. Following (19), we obtain

$$W_{k,Y_j/\tau_i} = \sum_{l=1}^{j} N_{kl} \sqrt{\frac{(\tau_i - Y_j)^2(Y_l - Y_{l-1})}{\tau_i(\tau_i - Y_l)(\tau_i - Y_{l-1})}}, \qquad Y_j \leq \tau_i,$$

$$W_{k+3,(Y_j-\tau_i)/(t-\tau_i)} = \sum_{l=1}^{j} N_{kl} \sqrt{\frac{(\tau_i - Y_j)^2(Y_l - Y_{l-1} \vee \tau_i)}{(t - \tau_i)(\tau_i - Y_l)(\tau_i - Y_{l-1} \vee \tau_i)}} \mathbb{I}[Y_l > \tau_i],$$

$$Y_j > \tau_i,$$

for $1 \leq j \leq \Lambda$, $k = 1,2,3$, and for both $i = 1,2$. Notice that there are alternative choices for the realization of the required locations of the standard BBs. However, the above formula has been carefully developed to ensure continuity in the parameter $\theta$. These expressions allow us now to write the proposed paths as deterministic functions of $E, \Psi, Z$ and $N$. In particular, substituting the above expressions into (18) and setting $\chi_{ij}$ equal to the right-hand side of (18) for $\tau = \tau_i$, with $i = 1, 2, 1 \leq j \leq \Lambda$, we obtain the final estimator:

$$a_i(\Xi, \theta) = p_i \times \prod_{j=1}^{\Lambda}[1 - \phi(\chi_{ij}, \theta)/\lambda], \qquad i = 1, 2.$$

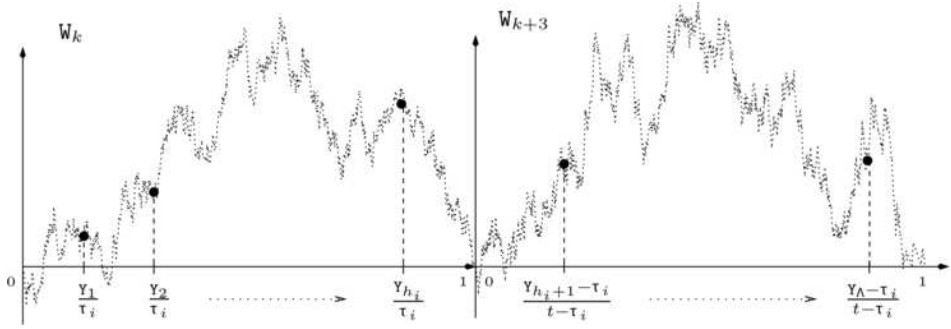

FIG. 1. *Scheme for simulating* $W_k$ *and* $W_{k+3}, k = 1,2,3$, *at the required discrete time instances; we have defined* $h_i = \sum_{j=1}^{\Lambda} \mathbb{I}[Y_j \leq \tau_i]$.

6.5. *A.s. continuity of the random function.* (Cnt2)(i) implies that $\beta_{ij}, \gamma_{ijl}$ are continuous functions of $\theta$ for all $i, j, l$ (even at values of $\theta$ such that $Y_j = \tau_i$ or $Y_l = \tau_i$); so, $\chi_{ij}$ is continuous in $\theta$. (Cnt1)(i), (Cnt2)(iii) suggest that $\phi(\cdot, \cdot)$ is continuous on $\mathbf{R} \times \Theta$, thus, $\phi(\chi_{ij}, \theta)$ and $a(\Xi, \theta)$ are continuous in $\theta$. The continuity of $L(\Xi, \theta)$ follows from assumptions (Cnt1)(ii) and (Cnt2).

6.6. *A.s. boundedness of the random function.* Note that $a(\Xi, \theta) \leq 1$ a.s., thus, w.p.1 $|L(\Xi, \theta)| < M$ for all $\theta \in \Theta$, where

$$M = \sup_{\theta \in \Theta}\{|\eta'(w,\theta)|\mathcal{N}_t(y-x)\exp\{A(y,\theta) - A(x,\theta) - l(\theta)t\}\} < \infty.$$

$M$ is finite as the supremum of a continuous function over a compact set.

**7. Discussion.** In this paper we have introduced a new method for likelihood inference for discretely observed diffusions. The method is computationally efficient and simple to implement, and expands significantly the family of diffusion models for which routine maximum likelihood calculations are possible. Applications of the approach advocated here together with other likelihood methods based on EA are given in [8].

Our methodology is a Monte Carlo approach based on two types of probabilistic constructions. The first of these exploits a duality between diffusions and Poisson processes, which has become transparent since the development of EA in [7]. The second involves explicit representations of conditioned Brownian sample paths constructed so as to be suitably continuous in model parameters. Although the entire mathematical construction is quite involved, in this paper we are able to distill its implementation down to a very simple Monte Carlo algorithm requiring merely collections of Gaussian and exponential random variables, as described in Theorem 1. The computer code for implementing our algorithm is freely available by e-mail request to any of the authors.

The scope of the methodology we introduce here can be extended in various directions. A direct extension is to certain multivariate diffusions. For example, it can be seen that our methods extend to processes which (possibly after a transformation) can be written as

$$dX_s = \nabla A(X_s; \theta)\,ds + dB_s, \qquad X_s \in \mathcal{X} \subseteq \mathbf{R}^{\mathbf{m}},$$

and $B$ is a $m$-dimensional Brownian motion. The conditions $(C_0)$–$(C_3)$ have to be appropriately modified, particularly $(C_3)$ now becomes that $\alpha = \nabla A$ is such that $(\alpha^2 + \alpha')(\cdot;\theta)$ is bounded above on $(z_1, \infty) \times (z_2, \infty) \times \cdots \times (z_m, \infty)$ for all $z_1, \ldots, z_m \in \mathbf{R}$. Furthermore, it is likely that our methods will be extended by currently ongoing work on simulating diffusions for which condition $(C_3)$ is not required to hold.

## APPENDIX: SIMULATION OF $\tau$

The joint distribution of the minimum and the time instance when this is attained $(m, \tau)$ for a Brownian bridge $W$ is of a known form; see [36]. In [7] we noticed that the distribution of $\tau$ conditionally on $m$ can be expressed in terms of the mixture of two inverse Gaussian laws. The exact simulation formula is the following:

$$\tau^{-1} = \{1 + \mathbb{I}[U_0 < (1 + \sqrt{c_1/c_2})^{-1}] \cdot I_1 + \mathbb{I}[U_0 \geq (1 + \sqrt{c_1/c_2})^{-1}] \cdot I_2^{-1}\}/t,$$

where $c_1 = (y-m)^2/(2t)$, $c_2 = (x-m)^2/(2t)$, and $I_1 \sim IGau(\sqrt{c_1/c_2}, 2c_1)$, $I_2 \sim IGau(\sqrt{c_2/c_1}, 2c_2)$ and $U_0 \sim \text{Un}[0,1]$. We denoted by $IGau(\cdot, \cdot)$ the inverse Gaussian distribution specified by two positive parameters (see, e.g., Chapter IV.4 of [14]). An inverse Gaussian random variable $I$ with parameters $(c, d)$ can be represented as (see page 149 of [14])

$$I = \mathbb{I}[V_0 \leq 1/(1+z)]cz + \mathbb{I}[V_0 > 1/(1+z)]\frac{c}{z},$$

where $z = 1 + (c/d)\mathtt{Z}^2/2 - \sqrt{4(c/d)\mathtt{Z}^2 + (c/d)^2\mathtt{Z}^4}/2$, $V_0 \sim \text{Un}[0,1]$. Using this formula for $I_1, I_2$ with the same $V_0, \mathtt{Z}$, and after some algebra, the expression for $\tau$ simplifies to (17).

**Acknowledgments.** We thank the referees for their valuable suggestions and corrections.

## REFERENCES

[1] Aït-Sahalia, Y. (2002). Maximum likelihood estimation of discretely sampled diffusions: A closed-form approximation approach. *Econometrica* **70** 223–262. MR1926260
[2] Andrews, D. W. K. (1987). Consistency in nonlinear econometric models: A generic uniform law of large numbers. *Econometrica* **55** 1465–1471. MR0923471
[3] Asmussen, S., Glynn, P. and Pitman, J. (1995). Discretization error in simulation of one-dimensional reflecting Brownian motion. *Ann. Appl. Probab.* **5** 875–896. MR1384357
[4] Attouch, H. (1984). *Variational Convergence for Functions and Operators*. Pitman (Advanced Publishing Program), Boston, MA. MR0773850
[5] Beck, A. (1963). On the strong law of large numbers. In *Ergodic Theory (Proc. Internat. Sympos., Tulane Univ., New Orleans, La., 1961)* 21–53. Academic Press, New York. MR0160256
[6] Beskos, A., Papaspiliopoulos, O. and Roberts, G. O. (2008). A factorisation of diffusion measure and finite sample path constructions. *Methodol. Comput. Appl. Probab.* **10** 85–104. MR2394037
[7] Beskos, A., Papaspiliopoulos, O. and Roberts, G. O. (2006). Retrospective exact simulation of diffusion sample paths with applications. *Bernoulli* **12** 1077–1098. MR2274855
[8] Beskos, A., Papaspiliopoulos, O., Roberts, G. O. and Fearnhead, P. (2006). Exact and computationally efficient likelihood-based estimation for discretely observed diffusion processes (with discussion). *J. R. Stat. Soc. Ser. B Stat. Methodol.* **68** 1–29. MR2278331
[9] Beskos, A. and Roberts, G. O. (2005). Exact simulation of diffusions. *Ann. Appl. Probab.* **15** 2422–2444. MR2187299

[10] BIBBY, B. M. and SØRENSEN, M. (1995). Martingale estimation functions for discretely observed diffusion processes. *Bernoulli* **1** 17–39. MR1354454

[11] BURKE, G. (1965). A uniform ergodic theorem. *Ann. Math. Statist.* **36** 1853–1858. MR0183001

[12] CAPPÉ, O., DOUC, R., MOULINES, E. and ROBERT, C. (2002). On the convergence of the Monte Carlo maximum likelihood method for latent variable models. *Scand. J. Statist.* **29** 615–635. MR1988415

[13] CHOIRAT, C., HESS, C. and SERI, R. (2003). A functional version of the Birkhoff ergodic theorem for a normal integrand: A variational approach. *Ann. Probab.* **31** 63–92. MR1959786

[14] DEVROYE, L. (1986). *Nonuniform Random Variate Generation*. Springer, New York. MR0836973

[15] DUDLEY, R. M. (1999). *Uniform Central Limit Theorems*. Cambridge Univ. Press. MR1720712

[16] DURHAM, G. B. and GALLANT, A. R. (2002). Numerical techniques for maximum likelihood estimation of continuous-time diffusion processes (with discussion). *J. Bus. Econom. Statist.* **20** 297–338. MR1939904

[17] ELERIAN, O., CHIB, S. and SHEPHARD, N. (2001). Likelihood inference for discretely observed nonlinear diffusions. *Econometrica* **69** 959–993. MR1839375

[18] ERAKER, B. (2001). MCMC analysis of diffusion models with application to finance. *J. Bus. Econom. Statist.* **19** 177–191. MR1939708

[19] FERGUSON, T. S. (1996). *A Course in Large Sample Theory*. Chapman and Hall, London. MR1699953

[20] GALLANT, A. R. and TAUCHEN, G. (1996). Which moments to match? *Econometric Theory* **12** 657–681. MR1422547

[21] GENON-CATALOT, V. and JACOD, J. (1993). On the estimation of the diffusion coefficient for multi-dimensional diffusion processes. *Ann. Inst. H. Poincaré Probab. Statist.* **29** 119–151. MR1204521

[22] GEYER, C. J. (1994). On the convergence of Monte Carlo maximum likelihood calculations. *J. Roy. Statist. Soc. Ser. B* **56** 261–274. MR1257812

[23] GEYER, C. J. and THOMPSON, E. A. (1992). Constrained Monte Carlo maximum likelihood for dependent data (with discussion). *J. Roy. Statist. Soc. Ser. B* **54** 657–699. MR1185217

[24] GIESY, D. P. (1976). Strong laws of large numbers for independent sequences of Banach space-valued random variables. In *Probability in Banach Spaces (Proc. First Internat. Conf., Oberwolfach, 1975)* 89–99. *Lecture Notes in Math.* **526**. Springer, Berlin. MR0440673

[25] GOURIEROUX, C., MONFORT, A. and RENAULT, E. (1993). Indirect inference. *J. Applied Econometrics* **8** 85–118.

[26] HALL, P. and HEYDE, C. C. (1980). *Martingale Limit Theory and Its Application*. Academic Press, New York. MR0624435

[27] KESSLER, M. and PAREDES, S. (2002). Computational aspects related to martingale estimating functions for a discretely observed diffusion. *Scand. J. Statist.* **29** 425–440. MR1925568

[28] KINGMAN, J. F. C. (1993). *Poisson Processes*. Oxford Univ. Press, New York. MR1207584

[29] KLOEDEN, P. E. and PLATEN, E. (1992). *Numerical Solution of Stochastic Differential Equations*. Springer, Berlin. MR1214374

[30] MOURIER, E. (1953). Eléments aléatoires dans un espace de Banach. *Ann. Inst. H. Poincaré* **13** 161–244. MR0064339

[31] ØKSENDAL, B. (2003). *Stochastic Differential Equations*, 6th ed. Springer, Berlin. MR2001996
[32] PEDERSEN, A. R. (1995). Consistency and asymptotic normality of an approximate maximum likelihood estimator for discretely observed diffusion processes. *Bernoulli* **1** 257–279. MR1363541
[33] PESKIR, G. (2000). *From Uniform Laws of Large Numbers to Uniform Ergodic Theorems*. Univ. Aarhus, Dept. Mathematics. MR1805157
[34] PRESS, W. H., TEUKOLSKY, S. A., VETTERLING, W. T. and FLANNERY, B. P. (1992). *Numerical Recipes in C*, 2nd ed. Cambridge Univ. Press. MR1201159
[35] ROBERTS, G. O. and STRAMER, O. (2001). On inference for partially observed nonlinear diffusion models using the Metropolis–Hastings algorithm. *Biometrika* **88** 603–621. MR1859397
[36] SHEPP, L. A. (1979). The joint density of the maximum and its location for a Wiener process with drift. *J. Appl. Probab.* **16** 423–427. MR0531776
[37] SØRENSEN, H. (2004). Parametric inference for diffusion processes observed at discrete points in time: A survey. *Internat. Statist. Rev.* **72** 337–354.
[38] TALAGRAND, M. (1987). The Glivenko–Cantelli problem. *Ann. Probab.* **15** 837–870. MR0893902
[39] TUCKER, H. G. (1959). A generalization of the Glivenko–Cantelli theorem. *Ann. Math. Statist.* **30** 828–830. MR0107891
[40] WALD, A. (1949). Note on the consistency of the maximum likelihood estimate. *Ann. Math. Statistics* **20** 595–601. MR0032169

A. BESKOS
O. PAPASPILIOPOULOS
G. ROBERTS
DEPARTMENT OF STATISTICS
UNIVERSITY OF WARWICK
COVENTRY CV4 7AL
UNITED KINGDOM
E-MAIL: a.beskos@warwick.ac.uk
o.papaspiliopoulos@warwick.ac.uk
gareth.o.roberts@warwick.ac.uk